\documentclass[a4paper,11pt]{amsart}
\usepackage{a4wide,amsfonts,amsmath,amssymb,amsthm,epsfig,enumerate,hyperref}
\usepackage{thm-restate}
\usepackage{color}
\usepackage{subcaption}
\usepackage{wrapfig}
\usepackage{listings}

\usepackage[textsize=footnotesize]{todonotes}

\graphicspath{{graphics/}}

\newcommand{\Cay}{\ensuremath{\mathrm{Cay}}}
\newcommand{\Aut}{\ensuremath{\mathrm{Aut}}}
\newcommand{\ops}{\ensuremath{\text{ops}}}
\newcommand{\FF}{\mathbb{F}}
\newcommand{\NN}{\mathbb{N}}
\newcommand{\ZZ}{\mathbb{Z}}
\newcommand{\cO}{\mathcal{O}}

\newcommand{\boxprod}{\mathbin{\square}}

\DeclareMathOperator{\col}{col}
\DeclareMathOperator{\poly}{poly}

\newcommand{\vspan}[1]{{\ensuremath{\left\langle #1 \right\rangle}}}

\definecolor{defblue}{rgb}{0.1,0.1,0.7}
\newcommand{\defi}[1]{\textcolor{defblue}{\emph{#1}}}

\newtheorem{thm}{Theorem}
\newtheorem{lem}[thm]{Lemma}

\newtheorem{cor}[thm]{Corollary}

\newtheorem{con}{Conjecture}

\theoremstyle{definition}

\newtheorem{defn}{Definition}

\newtheorem*{claim*}{Claim}
\newenvironment{claimproof}[1][\proofname]{
  
  \begin{proof}[#1]}{\end{proof}}

\title{Generating all invertible matrices by row operations}

\author[P.~Gregor]{Petr Gregor}
\address[P.~Gregor]{Department of Theoretical Computer Science and Mathematical Logic, Charles University, Prague, Czech Republic}
\email{gregor@ktiml.mff.cuni.cz}

\author[H.~P.~Hoang]{Hung. P. Hoang}
\address[H.~P.~Hoang]{Algorithm and Complexity Group, Faculty of Informatics, TU Wien, Austria}
\email{phoang@ac.tuwien.ac.at}

\author[A.~Merino]{Arturo Merino}
\address[A.~Merino]{Engineering Institute, Universidad de O'Higgins, Rancagua, Chile}
\email{arturo.merino@uoh.cl}

\author[O.~Mička]{Ondřej Mička}
\address[O.~Mička]{Department of Theoretical Computer Science and Mathematical Logic, Charles University, Prague, Czech Republic}
\email{micka@ktiml.mff.cuni.cz}

\thanks{This work was supported by Czech Science Foundation grant GA~22-15272S and the project TIPEA that has received funding from the European Research Council (ERC) under the European Unions Horizon 2020 research and innovation programme (grant agreement No. 850979).
Hung P. Hoang further acknowledges support from the Austrian Science Foundation (FWF, project Y1329 START-Programm).
This work was initiated at the 2nd Combinatorics, Algorithms, and Geometry workshop in Dresden, Germany in 2022.
We would like to thank the organizer and participants of the workshop for the inspiring atmosphere.}

\begin{document}
\begin{abstract}
  We show that all invertible $n \times n$ matrices over any finite field $\mathbb{F}_q$ can be generated in a \emph{Gray code} fashion. More specifically, there exists a listing such that (1) each matrix appears exactly once, and (2) two consecutive matrices differ by adding or subtracting one row from a previous or subsequent row, or by multiplying or dividing a row by the generator of the multiplicative group of $\mathbb{F}_q$.
  This even holds in the more general setting where the pairs of rows that can be added or subtracted are specified by an arbitrary transition tree that has to satisfy some mild constraints.
  Moreover, we can prescribe the first and the last matrix if $n\ge 3$, or $n=2$ and $q>2$.
  In other words, the corresponding flip graph on all invertible $n \times n$ matrices over $\mathbb{F}_q$ is Hamilton connected if it is not a cycle.
  This solves yet another special case of Lov\'{a}sz conjecture on Hamiltonicity of vertex-transitive graphs.

\end{abstract}

\maketitle

\section{Introduction}
Combinatorial generation is one of the most basic tasks we can perform on combinatorial objects and a key topic in Volume 4A of Knuth's seminal series \emph{The Art of Computer Programming}~\cite{MR3444818}.
In this task, we are given an implicit description of the objects and need to produce a listing of all objects fitting the description, with each object appearing exactly once.
The goal is to develop an algorithm that can generate these objects at a fast rate.

If consecutive objects produced by a generation algorithm differ by large changes, the algorithm must spend a lot of time updating its data structures.
Therefore, a natural first step towards creating an efficient generation algorithm is to ensure that consecutive objects differ by only a \emph{small change}.
Such a listing is known as a \defi{(combinatorial) Gray code}; see Mütze's survey~\cite{TheSurvey} for many Gray codes of various objects.
In addition to combinatorial generation, Gray codes are also relevant in the field of combinatorial reconfiguration, which examines the relationships between combinatorial objects through their local changes; see, e.g., Nishimura's recent introduction on  reconfiguration~\cite{MR4731617}.

In this paper, we study Gray codes for invertible matrices over a finite field.
A natural attempt for enumerating all invertible $n \times n$ matrices over a finite field $\mathbb{F}_q$, is to choose any nonzero first row and then selecting the following rows to be independent to the previous rows.
However, this attempt is not efficient as it requires multiple checks for independence to generate even a single matrix.
Furthermore, consecutive matrices in this listing may differ in multiple rows.
Instead, we focus on generating matrices in a Gray code order, i.e., every matrix is obtained from the previous one by a single elementary row operation.
We note that generating invertible matrices with specific properties has applications in cryptography, e.g., in McEliece cryptosystems \cite{MR3750441}.

\subsection{Strong Lov\'{a}sz conjecture}
All invertible $n \times n$ matrices over $\mathbb{F}_q$ with matrix multiplication form the \defi{general linear group} $GL(n,q)$.
Each elementary row operation can be represented by multiplying on the left by a matrix that corresponds to this row operation.
Hence, we are interested in finding a Hamilton path in an (undirected) Cayley graph on $GL(n,q)$ generated by the allowed row operations,
which is in turn an instance of Lovász conjecture~\cite{MR0263646} on the Hamiltonicity of vertex-transitive graphs.%
\footnote{A graph is \defi{vertex-transitive} if its automorphism group acts transitively on the vertices.}

Stronger versions of Lovász conjecture have been considered in the literature.
For example, Dupuis and Wagon \cite{Dupuis2015LaceableK} asked which non-bipartite vertex-transitive graphs are not Hamilton connected.
A graph is \defi{Hamilton connected} if there is a Hamilton path between any two vertices.
Similarly, they asked which bipartite vertex-transitive graphs are not Hamilton laceable \cite{Dupuis2015LaceableK}.
A bipartite graph is \defi{Hamilton laceable} if there is a Hamilton path between any two vertices from different bipartite sets.
Note that the bipartite sets must be of equal size, which is true for all vertex-transitive bipartite graphs except $K_1$.

\begin{con}[Strong Lov\'{a}sz conjecture]
\label{con:strongLovasz}
  For every finite connected vertex-transitive graph $G$ it holds that $G$ is Hamilton connected, or Hamilton laceable, or a cycle, or one of the five known counterexamples.
\end{con}

The five known counterexamples are the dodecahedron graph, the Petersen graph, the Coxeter graph, and the graphs obtained from the latter two by replacing each vertex with a triangle.
The dodecahedron graph is a non-bipartite vertex-transitive graph that has a Hamilton cycle, but it is not Hamilton connected~\cite{Dupuis2015LaceableK}.
The other four well-known counterexamples are non-bipartite vertex-transitive graphs that do not admit a Hamilton cycle.
Note that except when $G \in \{K_1, K_2, C_3, C_4\}$ the cases in the conjecture are mutually exclusive.

There are many results in line with Conjecture~\ref{con:strongLovasz}.
Particularly relevant to us is a result of Tchuente~\cite{MR683982} showing that the Cayley graph of the symmetric group $S_n$, generated by any connected set of transpositions, is Hamilton laceable when $n\geq 4$.
Another relevant example is Chen and Quimpo's Theorem~\cite{MR641233} showing that all Abelian Cayley graphs satisfy Conjecture~\ref{con:strongLovasz}.
Nevertheless, Conjecture~\ref{con:strongLovasz} remains open even for Cayley graphs of the symmetric group with every generator an involution~\cite{MR1201997}.
Note that none of the five counterexamples to Conjecture~\ref{con:strongLovasz} is a Cayley graph, leading to Cayley graph variants of Conjecture~\ref{con:strongLovasz} (e.g.,~\cite{MR1298976}).

\subsection{Row operations}
Our aim when generating all invertible matrices by row operations is to restrict the allowed operations as much as possible.
Note that for $q>2$ we must allow row multiplications by some scalar to be able to generate all $1 \times 1$ matrices.
Thus, we allow row multiplications by a fixed generator $\alpha$ of the multiplicative group of nonzero elements of $\mathbb{F}_q$.
We will also allow row multiplication by $\alpha^{-1}$; i.e., division by $\alpha$, to have an inverse operation for an undirected version of the problem.
Furthermore, we specify allowed row additions and subtractions by a directed \defi{transition graph} $T$ on the vertex set $[n]$, where $[n] := \{1, \dots, n\}$.
An edge $(i,j)\in E(T)$ specifies that we can add to or subtract from the $j$-th row the $i$-th row.
Then, each allowed row operation above corresponds to the left multiplication by a corresponding matrix from a set $\ops(T)$, formally defined by \eqref{eq:ops}.

Observe that to generate all invertible matrices by the allowed operations, the transition graph $T$ \emph{must} be strongly connected; see Lemma~\ref{lem:strongly_connected} below.
For our main result we require the following stronger condition.

\begin{defn}[Bypass transition graph]
	A transition graph~$T$ on the vertex set $[n]$ is a \defi{bypass transition graph} if either (i) $n=1$, or (ii) $n \geq 2$ and
	\begin{itemize}
		\item there exist an edge $(i, n)$ and an edge $(n, j)$ for some $i, j \in [n-1]$, and
		\item the graph $T - n$ obtained by removing $n$ from $T$ is also a bypass transition graph.
	\end{itemize}
\end{defn}

In other words, a bypass transition graph is obtained from a single vertex $1$ by repeatedly adding a directed path (a `\emph{bypass}') from some vertex $i$ to some vertex $j$ via a new vertex $n$.
An example of a transition graph with the above property is the one comprised by edges $(i, i+1)$ and $(i+1, i)$ for all $i\in [n-1]$; i.e., a bidirectional path.
In the language of row operations, this corresponds to allowing row additions or subtractions between any two consecutive rows.
It can be easily seen by induction that a bypass transition graph is strongly connected.

\subsection{Our results}
For any integer $n\ge 1$, a finite field $\mathbb{F}_q$, and an $n$-vertex transition graph $T$ we define the (undirected) Cayley graph
$$G(n,q,T):=\Cay(GL(n,q),\ops(T)),$$
where the set $\ops(T)$ is given by \eqref{eq:ops}.
Our main result is as follows.

\begin{thm}
\label{thm:ham_connected}
	Let $n\ge 2$ be an integer and $q$ be a prime power such that $q \geq 3$ if $n=2$.
    Let~$T$ be an $n$-vertex bypass transition graph.
	Then the graph $G(n,q,T)$ is Hamilton connected.
\end{thm}

Note that for $n=1$ the transition graph $T$ has no edges, so $G(1,q,T)$ for any $q\ge 3$ is simply a $(q-1)$-cycle and $G(1,2,T)=K_1$.
For $n=q=2$, we have that $T=(\{1,2\},\{(1,2),(2,1)\})$ is the only bypass transition graph, so $G(n,q,T)$ is a $6$-cycle, which is not Hamilton connected.
Thus, we may restate our result as follows.

\begin{cor}
\label{cor:strongLovasz}
	Let $n\ge 1$ be an integer and $q$ be a prime power, and let~$T$ be an $n$-vertex bypass transition graph.
	Then the graph $G(n,q,T)$ is Hamilton connected unless it is a cycle.
\end{cor}

\begin{figure}
	\centering
	\includegraphics[width=\textwidth]{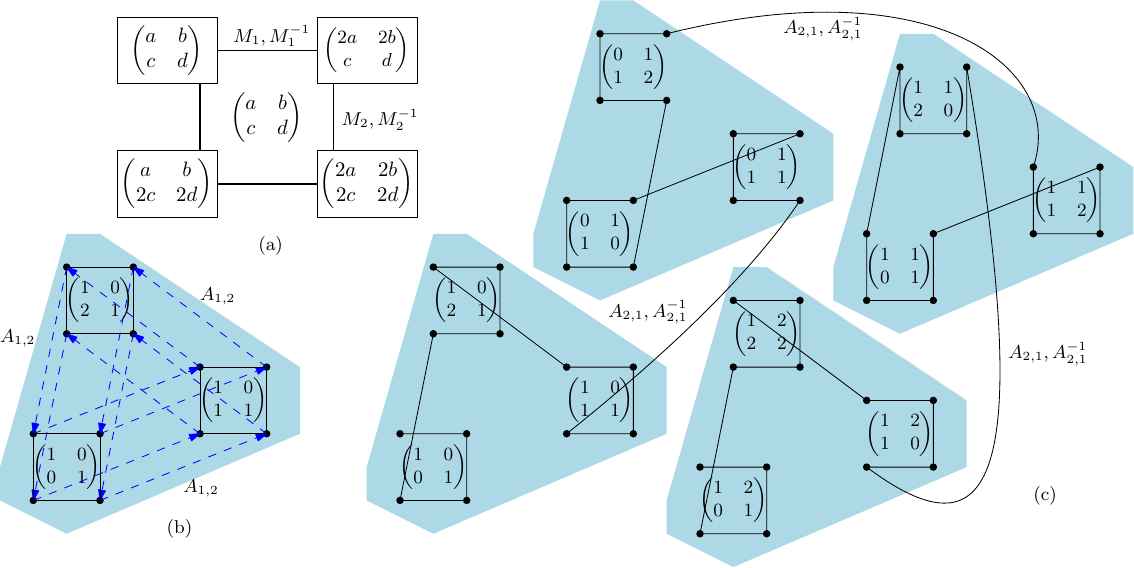}
	\caption{
		The part (c) illustrates a Hamilton path in the graph $G(2,3,([2],\{(1,2),(2,1)\}))$.
		Four vertices around a matrix~$Z$ are those obtained from $Z$ by multiplying or dividing a row by $\alpha$ (which is 2 for $q = 3$); see the part (a).
		The part (b) shows the edges within a shaded component, where the black solid edges are row multiplications/divisions, while the (directed) dashed edges are additions from the first row to the second row.
		Note that the other directions of the latter edges indicate subtractions of the first row from the second row.
		Furthermore, while these shaded components exhibit a Cartesian product structure, the same does not hold for the whole graph.
	}
	\label{fig:G23}
\end{figure}

This shows that the family of graphs $G(n,q,T)$ where $T$ is a bypass transition graph is yet another example of a family of Cayley graphs satisfying Conjecture~\ref{con:strongLovasz}.
A particularly interesting example is when $T$ is a bidirectional path.
See Figure~\ref{fig:G23} for an illustration for $n=2$ and $q=3$.

Moreover, we discuss how to turn the proof of Theorem~\ref{thm:ham_connected} algorithmic in Section~\ref{sec:algorithmic}.

\subsection{Related work}
Permutations of $[n]$ can be represented as (invertible binary) permutation matrices forming a subgroup of $GL(n,2)$.
Thus, all the vast results on generating permutations such as in \cite{MR4060409,MR683982} can be directly translated into the context of generating permutation matrices.
In particular, there is a general permutation framework developed in \cite{MR4391718} that allows us to generate many combinatorial classes by encoding them into permutations avoiding particular patterns.
However, the row operations that we consider here do not preserve the subgroup of permutation matrices, so our results do not fall into this framework.

A related task to generation is random sampling.
The construction of a random invertible $n \times n$ matrix over $\mathbb{F}_q$ is usually done by constructing a uniformly random matrix and checking whether it is non-singular.
The success probability is lower-bounded by a constant independent of $n$ but dependent on $q$ (e.g., see~\cite{Cooper} and the citations therein).
Hence, there is only a constant factor overhead for random sampling of an invertible matrix over a finite field compared to that of a matrix over the same field.
The latter task can be achieved, for example, by independently constructing each row (or column).

\section{Preliminaries}
\label{sec:prelim}
The \defi{(undirected) Cayley graph} of a group $\Gamma$ with a generator set $S$ is the graph $\Cay(\Gamma,S):=(\Gamma,\{\{x,sx\}: x\in\Gamma, s\in S\})$, assuming that $S$ is closed under inverses and does not contain the neutral element.
Note that we apply generators on the left as it is more natural for row operations on matrices.

The \defi{general linear group} $GL(n,q)$ is the group of all invertible $n \times n$ matrices over the finite field $\mathbb{F}_q$ with matrix multiplication.
Note that for $\mathbb{F}_q$ to be a field, $q$ has to be a prime power.
For example, $GL(1,2)$ is the trivial group, $GL(2,2)\simeq S_3$, and $GL(3,2)\simeq PSL(2,7)$ is also known as the group of automorphisms of the Fano plane. The number of elements in $GL(n,q)$ is $a_n:=(q^n-1)(q^n-q)\cdots(q^n-q^{n-1})$, which is obtained by counting choices for (nonzero) rows that are not spanned by the previous rows. It also satisfies the recurrence
\begin{equation}\label{eq:a_n}
a_{n}=(q^{n}-1)q^{n-1}a_{n-1},
\end{equation}
for $n\ge2$, and $a_1 = q - 1$ (i.e., the number of nonzero elements of $\mathbb{F}_q$).

By Gaussian elimination, the group $GL(n,q)$ can be generated by row additions and row multiplications by a scalar.
As we consider the Cayley graph to be undirected, we also consider the inverse operations, which we call row subtractions and row divisions by a scalar.
The formal definitions of these operations are as follows.

For $i\in [n]:=\{1,\dots,n\}$, let $r_i=r_i(A)$ denote the $i$-th row in $A$.
For distinct $x,y \in [n]$, we denote by $A_{xy}=(a_{ij})$ the binary matrix with $a_{ij}=1$ if and only if $i=j$, or ($i=y$ and $j=x$).
Note that left multiplication by $A_{xy}$ corresponds to adding the $x$-th row to the $y$-th row; i.e., the operation $r_x+r_y \to r_y$.
Similarly, multiplication by $A^{-1}_{xy}$ then corresponds to subtracting the $x$-th row to the $y$-th row; i.e., the operation $-r_x+r_y \to r_y$.

Let $\alpha$ be a generator of the multiplicative group of $\mathbb{F}_q$.
For $x \in [n]$, we denote by $M_{x} = (a_{ij})$ the matrix with $a_{ij} = \alpha$ if $i = j = x$, $a_{ij} = 1$ if $i = j \neq x$, and $a_{ij} = 0$ otherwise.
Left multiplication by $M_x$ corresponds to multiplying the $x$-th row by $\alpha$; i.e., the operation $\alpha r_x \to r_x$, and multiplication by $M^{-1}_{x}$ corresponds to the inverse operation $\alpha^{-1} r_x \to r_x$ that we call \defi{dividing} the $x$-th row by $\alpha$.
Note that for $q=2$ the multiplicative group is trivial, that is, $M_x=I$, where $I$ denotes the identity matrix.

A \defi{transition graph} $T$ is any directed graph on the vertex set $[n]$ with the edge set $E(T)$.
For a transition graph $T$ and a field $\mathbb{F}_q$ we define
\begin{equation}\label{eq:ops}
  \ops(T) := \{A_{ij}, A^{-1}_{ij} \, :\, (i,j) \in E(T) \} \cup \{M_i, M^{-1}_i \, :\, i \in [n]\},
\end{equation}
for $q>2$, and $\ops(T) := \{A_{ij}, A^{-1}_{ij} \, : \, (i,j) \in E(T) \}$ for $q=2$.
In other words, $\ops(T)$ contains the row additions and subtractions induced by the edges of $T$, and all row multiplications and divisions by $\alpha$ if they are nontrivial.
A directed graph is strongly connected if for any two vertices $i, j$, there is a directed path from $i$ to $j$.
A (strongly connected) component of a directed graph is a maximal induced subgraph that is strongly connected.
We make the following observation.
\begin{restatable}{lem}{stronglyconnected}
	\label{lem:strongly_connected}
	For every transition graph $T$, the set $\ops(T)$ generates the group $GL(n,q)$ if and only if $T$ is strongly connected.
\end{restatable}
\begin{proof}
If $T$ is not strongly connected, then there is a source component; i.e., there is no edge from a vertex outside the component to a vertex inside the component.
For a vertex $i$ in this source component, it is easy to see that the corresponding $i$-th rows of the matrices generated by the operations in $\ops(T)$ can only take value in the span of the rows indexed by this component.
Hence, $\ops(T)$ does not generate $GL(n,q)$.

Now assume that $T$ is strongly connected.
Observe that if we can add or subtract any row from any row and multiply or divide any row by $\alpha$, then by Gaussian elimination, we can generate any invertible matrix from any invertible matrix.
The row multiplications and divisions are already included in the definition of $\ops(T)$.
It remains to show that we can simulate any row addition or subtraction.
Suppose $A$ is the starting matrix.
For $a, b \in [n]$, $a \neq b$, by strong connectivity, there is a directed path $(v_1 = a, \dots, v_k = b)$ in $T$ for some $2 \leq k \leq n$.
We iteratively add the $v_i$-th row to the $v_{i+1}$-th row, for all $i \in [k-1]$.
At this point, the $v_i$-th row is equal to $\sum_{j \in [i]} r_{v_j}(A)$.
By repeatedly subtracting the $v_i$-th row from the $v_{i+1}$-th row, for $i = k-2, \dots, 1$, we restore the original values of the rows $r_{v_1}(A), \dots, r_{v_{k-1}}(A)$.
Next, we add the $v_i$-th row to the $v_{i+1}$-th row, for $i = 2, \dots, k-2$.
The $v_{k-1}$-th row is then $\sum^{k-1}_{i = 2} r_{v_i}(A)$.
Subtracting this row from the $k$-th row, the $k$-th row is then $r_1(A) + r_k(A)$.
Lastly, we subtract the $v_i$-th row from the $v_{i+1}$-th row, for $i = k-2, \dots, 2$.
The resulting matrix is equivalent to performing the operation $\{r_a + r_b \to r_b \}$ on $A$.
For the operation $\{- r_a + r_b \to r_b \}$, we perform exactly the same procedure as above, except for the very first operation, which is $\{ - r_{v_1} + r_{v_2} \to r_{v_2}\}$ instead of $\{r_{v_1} + r_{v_2} \to r_{v_2}\}$.
The lemma then follows.
\end{proof}

We denote by $\mathbb{F}_q^n$ the vector space of all $n$-tuples over the field $\mathbb{F}_q$.
The span of $u_1,\dots,u_k \in \mathbb{F}_q^n$ is denoted by $\langle u_1,\dots,u_k\rangle$. Its orthogonal space $\langle u_1,\dots,u_k\rangle^\bot$ is the kernel of the matrix with rows $u_1,\dots,u_k$.

For $k \geq 3$ we denote by $C_k$ a cycle on $k$ vertices, and for $k\in \{1,2\}$ we define $C_k$ as the complete graph $K_k$.
We also denote the path on $k$ vertices by $P_k$ for $k\ge 1$.
The \defi{Cartesian product} $G \mathbin{\square} H$ of two graphs $G$ and $H$ is the graph with the vertex set $V(G)\times V(H)$ and the edge set $\{(u,v)(u',v) \, : \, uu' \in E(G), v\in V(H)\} \cup \{(u,v)(u,v') \, : \, u \in V(G), vv'\in E(H)\}$.
For a graph $G$ and a subset $U$ of vertices, we denote by $G[U]$ the subgraph of $G$ induced by $U$.
Similarly, for a graph $G$ and two subsets of vertices $U_1,U_2\subseteq V$, we use $E[U_1,U_2]$ to denote the set of edges between $U_1$ and $U_2$, i.e., $E[U_1,U_2] = \{xy \in E : x \in U_1, y \in U_2 \}$.

For an edge-colored graph, a trail in a graph is \defi{alternating} if any two consecutive edges on the trail differ in color.

\section{Lemmas for Hamilton connectivity and laceability}
In this section, we present several useful lemmas for Hamilton connectivity or Hamilton laceability.
The first lemma on Cartesian product of cycles follows directly from a more general result on abelian Cayley graphs by Chen and Quimpo~\cite{MR641233}.
\begin{lem}
\label{lem:cartesian}
	For any $k \geq 2$ and $i_1, \dots, i_k\ge 2$ the graph $C_{i_1} \mathbin{\square} \dots \mathbin{\square} C_{i_k}$ is Hamilton connected if some $i_j$ is odd, and Hamilton laceable otherwise.
\end{lem}

We will also need a similar result for a Cartesian product of a path and even cycle. It is likely to be known, but we provide a proof for completeness.
\begin{lem}
\label{lem:pathxcycle}
For any $i\ge 2$ and $j\ge 2$ the graph $P_i \mathbin{\square} C_{2j}$ is Hamilton laceable.
\end{lem}
\begin{proof}
	The case $i=2$ is already covered by Lemma~\ref{lem:cartesian}.
	Now assume $i\ge 3$.
	Define $G := P_i \mathbin{\square} C_{2j}$.
	Let $G_k$ for $1\le k \le i$ denote the $k$-th copy of $C_{2j}$ in $G$.
	Let $x,y$ be two vertices of $G$ to be connected by a Hamilton path, assuming that $x\in V(G_k)$ and $y\in V(G_{k'})$ for some $1\le k \le k' \le i$.
	Since $G$ is bipartite and has an even number of vertices, this is only possible if $x$ and $y$ belong to different parts of $G$.

	First consider the case when $k = 1$ and $k' = i$.
	Choose a neighbor $y'$ of $y$ in $G_{k'}$, and let $P$ be a Hamilton path of $G_{k'}$ with endpoints $y$ and $y'$.
	Let $y''$ be the neighbor of $y'$ in $G_{k'-1}$.
	Note that $y''$ and $y$ are in the same part of $G$.
	By the inductive hypothesis, there exists a Hamilton $xy''$-path~$P'$ in the graph obtained from $G$ by removing $G_i$.
	Then concatenating $P$, $y'y''$, and $P'$ yields a Hamilton $xy$-path of $G$, as desired.

	For the remaining case, we can assume $k' < i$.
	Otherwise, this can still hold, after we reverse the indices of the copies of $C_{2j}$ in $G$ and swap the labels of $x$ and $y$.
	Then by the inductive hypothesis, there exists a Hamilton $xy$-path~$P'$ in the subgraph of $G$ obtained by removing $G_i$.
	In this subgraph the vertices of $G_{i-1}$ have degree three.
	Since $j \ge 2$, this implies the existence of an edge $ab$ in $G_{i-1}$ on the path $P$.
	Let $a'$ and $b'$ be the neighbors of $a$ and $b$, respectively, in $G_{i-1}$.
	Then $a'b'$ is an edge of $G_i$, and hence, there is a Hamilton $a'b'$-path~$P'$ of $G_i$.
	Replacing the edge $ab$ on $P$ with the edge $aa'$, the path $P'$, and the edge $b'b$ yields a Hamilton $xy$-path of $G$.
\end{proof}

The third lemma states that Hamilton connectivity of a nontrivial graph is preserved by a Cartesian product with any cycle. 
Note that it does not hold if $G=K_2$.
\begin{lem}
\label{lem:cycleproduct}
Let $G$ be a Hamilton connected graph on at least $3$ vertices. 
Then $G \mathbin{\square} C_k$ is Hamilton connected for any $k\ge 1$.
\end{lem}

\begin{proof}
  Let $G_i$ for $1\le i \le k$ denote the $i$-th copy of $G$ in $G \mathbin{\square} C_k$. Let $x,y$ be two vertices of $G \mathbin{\square} C_k$ to be connected by a Hamilton path, assuming that $x\in V(G_1)$ and $y\in V(G_j)$ for some $1\le j \le k$.

  Firstly, we connect all vertices of copies $G_i$ for $i=1, \dots, j$ into an $xy$-path. For this, we select vertices $x_i, y_i$ in $G_i$ such that $x_1=x$, $y_j=y$, $x_i\ne y_i$ for every $i=1,\dots, j$, and $y_i$ is a neighbor of $x_{i+1}$ for every $i=1,\dots, j-1$. Such vertices exist since $|V(G)|\ge 3$. Then we concatenate Hamilton paths in each $G_i$ between $x_i$ and $y_i$ that exist by Hamilton connectivity of $G$ into a single path $P$ between $x$ and $y$.

  Secondly, we extend the path $P$ to all vertices of copies $G_i$ iteratively for $i=j+1, \dots, k$. Let $ab$ be an edge of $G_{i-1}$ that belongs to the current path $P$, and let $a',b'$ be the neighbors of $a,b$ in $G_i$. By replacing the edge $ab$ with the edge $aa'$, a Hamilton $a'b'$-path in $G_i$, and the edge $b'b$, we extend the path $P$ to $G_i$. After the last step for $i=k$ we obtain a Hamilton $xy$-path.
\end{proof}

The last lemma joins many Hamilton connected graphs into a larger one.
This lemma seems quite versatile.
Not only is it useful in our proof in the next section, but it also allows us to easily reprove several classical results on Hamilton connectivity, for example for the permutahedron \cite{MR683982}.

\begin{lem}[Joining lemma]
\label{lem:join}
	Let $G$ be a graph with the vertex set partitioned into $k\ge 2$ disjoint subsets $V_1, \dots, V_k$ such that following conditions hold.
	\begin{enumerate}
		\item[(1)] $G[V_i]$ is Hamilton connected for every $i \in [k]$;
		\item[(2)] Every vertex in every set $V_i$ has a neighbor in some different set $V_j$;
		\item[(3)] There are at least three pairwise disjoint edges between every two sets $V_i$, $V_j$.\footnote{We could weaken the condition~(3) for $k\ge 4$ so that we only need two disjoint edges between all pairs of sets except for two disjoint pairs.}

	\end{enumerate}
	Then $G$ is Hamilton connected.
\end{lem}

\begin{proof}
Let $x, y\in V$ be two vertices to be connected by a Hamilton path.
First we consider the case when they are in different sets $V_i$.
We can assume that $x \in V_1$ and $y\in V_k$, otherwise we rename the sets.
We select vertices $x_i, y_i\in V_i$ for every $i\in [k]$ so that $x_1=x$, $y_k=y$, $x_i\ne y_i$ for every $i\in [k]$, and $y_i$ is a neighbor of $x_{i+1}$ for every $i\in [k-1]$.
Such vertices exist since there are at least three edges between $V_i$ and $V_{i+1}$ for every $i\in [k-1]$ by the condition~(3).
Then we concatenate Hamilton paths in $G[V_i]$ between $x_i$ and $y_i$ for each $i=1,\dots,k$ that exist by the condition~(1) into a Hamilton $xy$-path in $G$. Note that in this case we did not use the condition $(2)$.

In the second case $x$ and $y$ are in the same set $V_i$.
We can assume that $x,y\in V_1$. Let $P$ be a Hamilton path in $G[V_1]$ between $x$ and $y$.
If $k=2$, let $ab$ be an edge of $P$ such that the neighbors $a'$ and $b'$ of $a$ and $b$ in $V_2$, respectively, are distinct.
Such neighbors exist by the condition~(2), and such an edge~$ab$ exists, because otherwise all vertices in $V_1$ are only adjacent to one vertex in $V_2$, a contradiction to the condition~(3).
By replacing the edge $ab$ with the edge $aa'$, a Hamilton path of $G[V_2]$ between $a'$ and $b'$, and the edge $b'b$ we obtain a Hamilton $xy$-path in $G$.

If $k>2$, let $ab$ be an edge of $P$ such that $a$ and $b$ have neighbors $a'$ and $b'$, respectively, in different sets $V_i$ for $i>1$.
Such an edge $ab$ exists since every vertex of $V_1$ has a neighbor in some other set $V_i$ by the condition~(2), and they cannot be all from the same set, say $V_2$, for otherwise, the condition~(3) for the sets $V_1$ and $V_3$ would not hold.
By the same argument as in the first case, there exists a Hamilton path $R$ between $a'$ and $b'$ in the subgraph $G[V_2\cup \cdots \cup V_k]$.
Finally, replacing the edge $ab$ on $P$ with the edge $aa'$, the path $R$, and the edge $b'b$ yields a Hamilton $xy$-path in $G$.
\end{proof}

\section{Proof of Theorem~\ref{thm:ham_connected}}
\label{sec:main_proof}

We will prove the theorem by induction on $n$, and let $G_n := G(n,q,T)$.
We say that an edge $\{X,A_{ij}X\}$ of $G_n$ is \defi{labeled} $ij$ and an edge $\{X,M_{i}X\}$ of $G_n$ is \defi{labeled} $i$.

The proof of the base cases for $q=2$ and $n=3$, and for $q\ge 3$ and $n=2$ is deferred to Section~\ref{sec:base}; see Lemmas~\ref{lem:base_case32} and~\ref{lem:base_case}.
Here, we prove the inductive step, so we assume that $q=2$ and $n\ge 4$, or $q\ge 3$ and $n\ge 3$, and that the statement holds for the graph $G(n-1,q,T-n)$.
Our main tool is the joining lemma from the previous section (Lemma~\ref{lem:join}).

\begin{proof}[Proof of the inductive step]
We view rows of an invertible $n\times n$ matrix $A$ as an ordered basis $(r_1,\dots,r_n)$ of the vector space $\mathbb{F}_q^n$. The first $n-1$ rows span a subspace of dimension $n-1$ which is orthogonal to some subspace of dimension $1$. That is,
$$\langle r_1, \dots, r_{n-1}\rangle=\langle u \rangle^\bot$$
for a nonzero $u\in \mathbb{F}_q^n$ that satisfies $Ru^T = \textbf{0}$, where $R$ is the $(n-1) \times n$ matrix whose rows are $r_1, \dots, r_{n-1}$.

We denote by $S_u$ the set of all the $(n-1) \times n$ matrices whose rows form a basis of $\langle u \rangle^\bot$.
Observe that this operation gives a bijection between $S_u$ and $GL(n-1,q)$: Remove the $i$-th column of every matrix in $S_u$, where $i$ is the index of the first nonzero element of $u$ (which exists, as $u$ is not the zero vector).
Furthermore, for every matrix~$X$ in $S_u$, we can add any vector as the last row to form an $n \times n$ invertible matrix, as long as this added vector is independent of the rows of~$X$ (i.e., any vector in $\mathbb{F}_q^n \setminus \langle u \rangle^{\bot}$).

Note that this tallies with the count in~\eqref{eq:a_n}.
Recall that $a_{n-1}$ denotes the number of elements in $GL(n-1,q)$.
There are $(q^n -1) / (q-1)$ choices for the one-dimensional subspace $\langle u \rangle$.
For each subspace (with a representative basis $u$), $S_u$ has $a_{n-1}$ elements, due to the aforementioned bijection.
Lastly, for each matrix~$X$ in $S_u$, there are $q^{n-1}(q-1)$ possible last rows, which can be obtained by adding a linear combination of the rows of~$X$ to an initial last row and then multiplying the sum by a power of $\alpha$.
Together, we recover the recurrence statement~(\ref{eq:a_n}).

Following the above analysis, we prove the inductive step in four smaller steps.

First, given a one dimensional subspace with a basis $u$ and a vector $v$ independent of the rows of any matrix in $S_u$, we denote by the tuple $(S_u, v)$ the set of all matrices in $GL(n,q)$ formed by adding $v$ as the last row to each of the matrices in $S_u$.
Since the aforementioned bijection is an isomorphism of $G_n[(S_u, v)]$ and $G(n-1,q,T-n)$, by inductive hypothesis we conclude that $G_n[(S_u, v)]$ is Hamilton connected.

Before we proceed, we note that row multiplications, divisions, and row additions that do not involve the last row only transform a matrix into another matrix in the same set $(S_u, v)$ for some $u, v$.
Next, adding a row to the last row transforms a matrix in $(S_u, v)$ into another matrix in $(S_u, v')$ for $v \ne v'$.
Lastly, adding the last row to another row transforms a matrix in $(S_u, v)$ into another matrix in $(S_{u'}, v)$, where $\langle u \rangle \ne \langle u' \rangle$.

Second, we denote by $(S_u,\langle v \rangle)$ the set of all matrices in $GL(n,q)$ formed by adding any multiple of $v$ as the last row to each of the matrices in $S_u$. Since multiplication by $\alpha$ generates all nonzero elements of $\mathbb{F}_q$, the edges of label $n$ that multiply the last row form a cycle of length $q-1$.
Hence, the graph $G_n[(S_u,\langle v \rangle)] \simeq G_n[(S_u,v)] \mathbin{\square} C_{q-1}$, which is Hamilton connected by Lemma~\ref{lem:pathxcycle}.

Third, given a one dimensional subspace with a basis $u$, we denote by $(S_u, *)$ the set of all matrices in $GL(n,q)$ whose first $n-1$ rows form a matrix in $S_u$.
Here, we use Lemma~\ref{lem:join} to join the subgraphs $G_n[(S_u, \langle v \rangle)]$ for all applicable $v$ to prove the Hamilton connectivity of $G_n[(S_u, *)]$.
The joining edges between these subgraphs have label $in$ for $i \in [n-1]$ (such an $i$ is guaranteed by the bypass property of $T$).
In order to use the lemma, we show that all of its conditions hold.
The condition (1) follows the second step above.
The condition (2) is satisfied, because for every $u, v \in \mathbb{F}^n_q$ and a matrix~$X$ in $(S_u,v)\subseteq (S_u,\langle v \rangle)$, $A_{in}X$ is a neighbor of $X$ in $G_n$ and in $(S_u, \langle v + r_i(X)\rangle)$, a different set than $(S_u, \langle v\rangle)$.
For the condition (3), given $u$ and two distinct $v, v' \notin \langle u \rangle^\bot$ such that $\langle v \rangle \ne \langle v' \rangle$, we have that $v$ is a linear combination of a basis of $\langle u \rangle^\bot$ and $v'$, and consequently $v=x+av'$ for some nonzero $x\in \langle u\rangle^{\bot}$ and a nonzero $a\in \mathbb{F}_q$.
As $S_u$ contains all matrices whose rows form a basis in $\langle u \rangle^\bot$, there exist three matrices~$X_1,X_2$, and $X_3$ in $S_u$ such that their $i$-th row is $x=v - av'$.
We can guarantee three matrices in $S_u$, because in the inductive step, $n \geq 3$ and $q\ge 3$, or $n\geq 4$ and $q=2$, and hence, when we fix the $n-2$ rows including the $i$-th row, there are at least three different choices for the remaining row.
Then the edges $\{ X_1, A_{in} X_1 \}$, $\{ X_2, A_{in} X_2 \}$, and $\{ X_3, A_{in} X_3 \}$ are the three distinct edges as required by the condition (3).
We can now apply Lemma~\ref{lem:join} and conclude that $G_n[(S_u, *)]$ is Hamilton connected.

Last, we again apply Lemma~\ref{lem:join} to join the different subgraphs $G_n[(S_u, *)]$ for all subspaces $\langle u \rangle$ to complete the inductive step.
Here, the joining edges have the label $nj$ for some $j \in [n-1]$, which exist because $T$ is a bypass transition graph.
The condition (1) of the lemma follows the previous step.
The condition (2) is satisfied, because for any~$X$ in some $(S_u, *)$, $A_{nj} X$ is a neighbor of $X$ in $G_n$ and belongs to a different set $(S_{u'}, *)$.
For the condition (3), given $u, u'$ not in the same one-dimensional subspace, $\langle u, u' \rangle^\bot$ is a subspace of dimension $n-2$.

If $n \geq 4$, or $n=3$ and $q>3$, there exist three distinct matrices $B$, $B'$, and $B''$ whose rows form bases of this $(n-2)$-dimensional subspace.
The remaining case $n=3$ and $q=3$ is considered separately below.
Let $v_u \in \langle u \rangle^\bot \setminus \langle u' \rangle^\bot$ and $v_{u'} \in \langle u' \rangle^\bot \setminus \langle u \rangle^\bot$.
Clearly, we have that $v_{u'} - v_u$ is independent of the rows of each matrix $B$, $B'$, and $B''$.
Let $\tilde{B}$, $\tilde{B}'$, and $\tilde{B}''$ be the $n \times n$ matrix obtained from $B$, $B'$, and $B''$ respectively by inserting a new row $v_u$ at the $j$-th position and $v_{u'} - v_u$ as the last row.

If $n=3$ and $q=3$ we have $\langle u \rangle^\bot \cap \langle u' \rangle^\bot=\langle w \rangle = \{0,w,2w\}$ for some nonzero $w\in \mathbb{F}_3^3$, so there are only two distinct matrices whose rows form bases of this $1$-dimensional subspace, in particular $B=(w)$ and $B'=(2w)$.
In this case, we define $\tilde{B}$ and $\tilde{B}'$ as in the previous case, but for $\tilde{B}''$ we take the matrix obtained from $B$ by inserting a new row $2v_u$ at the $j$-th position and $2v_{u'}-2v_u$ as the last row.

In both cases, $\{\tilde{B}, A_{nj} \tilde{B}\}$, $\{\tilde{B}', A_{nj} \tilde{B}'\}$, and $\{\tilde{B}'', A_{nj} \tilde{B}''\}$ are the three edges as required by the condition (3).
This completes all the conditions of Lemma~\ref{lem:join} and completes the inductive step.
\end{proof}

\section{Base cases}
\label{sec:base}

In this section, we prove the base cases for our induction.
First, we prove the following with computer assistance.

\begin{lem}
\label{lem:base_case32}
	For every bypass transition graph $T$, the graph $G(3,2,T)$ is Hamilton connected.
\end{lem}

\begin{proof}
	We verified the statement using computer search in SageMath; see Appendix~\ref{app:sage}.
\end{proof}

Since the only bypass transition graph $T$ for $n=2$ is the complete graph $T=([2],\{(1,2),(2,1)\})$, the remaining base cases for $q \geq 3$ and $n =2 $ are captured in the following lemma.
\begin{lem}
	\label{lem:base_case}
	For a prime power $q \geq 3$, $G(2,q,([2],\{(1,2),(2,1)\}))$ is Hamilton connected.
\end{lem}

To prove Lemma~\ref{lem:base_case}, we first consider the graph that arises by removing the edges that add the second row to first row; i.e., for a prime power $q \in \NN$, we define $G'(q) := G(2,q,([2],\{(1,2)\}))$.
The graph $G'(q)$ is disconnected, and thus, we consider the connected component in $G'(q)$ that contains the identity matrix and denote it by $H(q)$.
We will usually write $H$ and $G'$ for $H(q)$ and $G'(q)$ whenever there is no risk of confusion.

It is easy to see that the vertices of $H$ are of the following form:
\[
V(H) = \left\{ \begin{pmatrix} \alpha^i & 0 \\ \alpha^j a  & \alpha^j\end{pmatrix} : i,j \in \{0,\dots, q-2\}, a\in \FF_q \right\}.
\]
Furthermore, the graph $H$ has a simple structure when analyzing the components by fixing $a \in \FF_q$, as described in the following.
Let us define
\[
V_a = \left\{ \begin{pmatrix} \alpha^i & 0 \\ \alpha^j a  & \alpha^j\end{pmatrix} : i,j \in \{0,\dots, q-2\} \right\}
\]
and let $H_a$ be the graph induced by fixing $a$ in $H$; i.e.,
$H_a := H [V_a]$.
We have the following simple observations regarding $H$ and its decomposition by fixing $a\in \FF_q$.

\begin{enumerate}
	\item[(p1)] \emph{($H$ splits into copies of $H_a$.)}
	Removing the edges that add the first row to the second row in $H$ disconnects $H$ and splits it into the connected components $\{H_a :a \in \FF_q\}$.
	\item[(p2)] \emph{(The graphs $H_a$ have good Hamiltonicity properties.)}
	For every $a\in \FF_q$, we get that $H_a$ is a toroidal grid of dimensions $(q-1) \times (q-1)$ where each dimension of the grid is given by multiplication by $\alpha$ in the respective row; i.e., $H_a \cong C_{q-1} \boxprod C_{q-1}$.
	In particular, $H_a$ is isomorphic to $H_b$ for every $a,b\in \FF_q$.
	\item[(p3)] \emph{(The components $H_a$ are well-connected.)}
	For every $i \in \{0,\dots, q-2\}$ and $a,b \in \FF_q$ such that $a\neq b$,  we have that
	\begin{enumerate}
	\item $\alpha^i\begin{pmatrix} a-b & 0 \\
		a & 1
	\end{pmatrix} \in V_a$ and $\alpha^i\begin{pmatrix} a-b & 0 \\
		b & 1
	\end{pmatrix} \in V_b$ are connected by an edge,
	\item $\alpha^i\begin{pmatrix} b-a & 0 \\
		a & 1
	\end{pmatrix} \in V_a$ and $\alpha^i\begin{pmatrix} b-a & 0 \\
		b & 1
	\end{pmatrix} \in V_b$ are connected by an edge,
	\end{enumerate}
	and no other edges between~$V_a$ and~$V_b$ exist.
	In particular, for every  $a,b \in \FF_q$ such that $a\neq b$ we have that $|E[V_a,V_b]| = 2(q-1)$ with all the edges being disjoint.
\end{enumerate}

We exploit these properties as follows:
We split $H$ into the components $H_a$ for $a \in \FF_q$, then since the graphs $H_a$ are either Hamilton laceable or connected and the components $H_a$ are well-connected we can glue the corresponding Hamilton paths in each $H_a$ to form a Hamilton path in $H$.

If $q$ is even, for every $a\in \FF_q$ the graphs $H_a$ are Hamilton connected.
This makes it easier to lift the Hamilton paths from $H_a$ to a Hamilton path in $H$.
However, when $q$ is odd, the picture is much different.
In particular, there are parity constraints given by the fact that for every $a\in \FF_q$ the graph $H_a$ is now bipartite.
To make this formal, we partition $V(H)$ into two colors.
We say that $x = \begin{pmatrix}
	\alpha^i & 0 \\ \alpha^j a & \alpha^j
\end{pmatrix}$ is \defi{blue} whenever $i+j$ is even, and it is \defi{red} whenever $i+j$ is odd.
We denote the \defi{color of a vertex} $x\in V(H)$ by $\col(x)$.
It is easy to show that if $x \in V_a$, $y \in V_b$ for some $a\neq b$ and there is an edge $xy \in E[V_a,V_b]$, then $\col(x)=\col(y)$.
Thus, for every edge $xy \in E[V_a,V_b]$ we can define its \defi{(edge) color} as $\col(xy):=\col(x)=\col(y)$.

If $q \equiv 3 \pmod{4}$, a simple computation shows that whenever there is a red or blue edge between $V_a$ and $V_b$ for $a,b \in \FF_q$ we also have an edge of the opposite color.
Thus, the coloring does not impose any extra restrictions.

The problematic case occurs whenever $q\equiv 1 \pmod{4}$.
In this case, all the edges between $V_a$ and $V_b$ have the same color for $a, b \in \FF_q$.
Hence, it is natural to consider the graph where we contract every $V_a$ for $a \in \FF_q$.
Thus, we obtain a new graph \defi{$\bar{K}_{q}$} with $\FF_q$ as vertices, and for the edges $xy \in E[V_a,V_b]$ we put a new edge $ab$ colored with $\col(xy)$.
This graph is a complete graph on $\FF_q$ where the coloring of the edges can be succinctly described as follows:
\begin{itemize}
	\item[(*)] For $x$ and $y$ in $\mathbb{F}_q$, the edge $xy$ has color red (blue), if there exists an odd (even)~$z \in \mathbb{Z}$, such that $x - y = \alpha^z$. (See Figure~\ref{fig:coloring-kq} for an example.)
\end{itemize}

{\makeatletter
\let\par\@@par
\par\parshape0
\everypar{}\begin{wrapfigure}[16]{r}{.3\textwidth}
	\centering
	\vspace{-2.6em}
	\includegraphics[width=.2\textwidth]{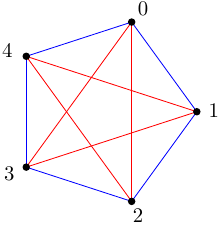}
	\caption{The graph $\bar{K}_5$ for $\alpha = 2$.}
	\label{fig:coloring-kq}
\end{wrapfigure}
\noindent
If we plan to have a Hamilton path of $H$ that traverses each set $V_a$ at a time for $a \in \FF_q$, then for any $a,b
\in \FF_q$, there is at most one edge of the Hamilton path that crosses between $V_a$ and $V_b$.
Further, as each $V_a$ has even size, this means that as we traverse this Hamilton path, any two consecutive such `crossing' edges have to differ in color.
This translates to the requirement that we should have an alternating Hamilton path in $\bar{K}_{q}$.
We show in the next lemma that this holds, even for Hamilton connectivity.

}%

\begin{restatable}{lem}{alternatingcoloring}
	\label{lem:alternating_coloring}
		Let $q$ be a prime power, $q \equiv 1 \pmod{4}$.
		For any two distinct vertices $a, b \in \mathbb{F}_q$ and a color $c$ of either red or blue, there exists an alternating Hamilton $ab$-path of $\bar{K}_{q}$, such that $a$ is incident to an edge of color~$c$ on the path.
\end{restatable}
We defer the proof of Lemma~\ref{lem:alternating_coloring} to Section~\ref{sec:alternating_coloring}.
We can use Lemma~\ref{lem:alternating_coloring} to prove Hamiltonicity properties of subgraphs of $H$.
To this end, we have the following definition.

\begin{defn}
\label{def:structured}
An induced subgraph $H'$ of $H$ is \defi{structured} if and only if the following holds:
\begin{enumerate}
	\item For every $a \in \FF_q$ we have that $H'[V_a]$ is isomorphic to either $C_{q-1} \boxprod C_{q-1}$ or $C_{q-1 } \boxprod P_{\ell}$ for $\ell \geq (q-1)/2$, and
	\item For every distinct $a, b \in \FF_q$, there is at least one edge between $H'[V_a]$ and $H'[V_b]$.
\end{enumerate}
\end{defn}

\begin{restatable}{lem}{ordercomponents}
	\label{lem:order_components}
	Let $q\geq 5$ be an odd integer and $H'$ be a structured induced subgraph of $H(q)$.
	Let $x, y\in V(H(q))$ be two vertices of different colors such that $x\in V_a$, $y\in V_b$ with distinct $a,b \in \FF_q$.
	Then there exists a Hamilton $xy$-path in $H'$.
\end{restatable}

The reader may notice similarities between Lemma~\ref{lem:order_components} and the joining lemma (Lemma~\ref{lem:join}).
In particular, they may wonder why we require only \emph{one} edge between components, instead of the three needed in the joining lemma.
Recall that the need for three edges in the joining lemma was in the case where we want to have an $xy$-subpath that spans two consecutive components, but the edges that cross between these two components are incident to either $x$ or $y$. 
However, this cannot happen for structured graphs, because the coloring conditions force these endpoints not to be used in crossing edges.%
	\begin{proof}[Proof of Lemma~\ref{lem:order_components}]
		For every $a \in \FF_q$, we define $H'_a := H'[V_a]$ and $V'_a := V(H'_a)$.
		Applying Lemma~\ref{lem:pathxcycle} and Lemma~\ref{lem:cycleproduct} we obtain that for every $a\in \FF_q$ the graph $H'_a$ is Hamilton laceable.
	
		Let $q' := (q-1)/2$ and assume without loss of generality that $x$ is red and $y$ is blue. 
		We consider two cases.
		
		\noindent\textbf{Case 1:} $q \equiv 3 \pmod{4}$.
		Suppose $c - d = \alpha^z$ for some~$z \in \{0, \dots, q-2\}$.
		Since $\alpha^{q'} = -1$, $d - c = \alpha^{z + q'}$.

		We first argue that there exist edges of both colors between $V'_c$ and $V'_d$, too.
		For $i,j \in \{0, \dots, q-2\}$ and $p \in \FF_q$, define $M(i,j,p) := \begin{pmatrix} \alpha^i & 0 \\ \alpha^j p & \alpha^j \end{pmatrix}$.
		Then the two types of edges between of $V_c$ and $V_d$ as described by (a) and (b) in (p3) are of the form $M(i+z, i, c) M(i+z, i, d)$ and $M(i+z+q', i, c) M(i+z+q', i, d)$ for some $i$.
		Further, since $q \equiv 3 \pmod{4}$, $q'$ is odd, and hence these two types of edges have different colors.
		Without loss of generality, suppose the edge guaranteed by condition (2) of Definition~\ref{def:structured} is of the first type, i.e., $M(i+z,i,c)M(i+z,i,d)$ for some $i \in \{0, \dots, q-2\}$.
		Then condition (1) of the Definition~\ref{def:structured} implies that there exist two subintervals $[t_1, t'_1]$ and $[t_2, t'_2]$ of the (cyclic) interval $[0,q-1]$, such that these subintervals have length at least $(q-1)/2$, and $M(j_1, j_2, c) \in V'_c$ and $M(j_1, j_2,d) \in V'_d$, for each $j_1 \in [t_1,t'_1]$ and $j_2 \in [t_2,t'_2]$.
		Since $[t_1,t'_1]$ and $[t_2+z+q',t'_2+z+q']$ have length at least $(q-1)/2$, they have at least one common integer point $j$.
		Then $M(j+z+q', j, c) M(j+z+q', j, d)$ is an edge between $V'_c$ and $V'_d$.
		As argued before, this edge has different color than that of $M(i+z,i,c)M(i+z,i,d)$.
		Overall, we have two edges of two colors between $V'_c$ and $V'_d$, as required.

		We now choose any permutation $a_1\cdots a_q$ of $\FF_q$ with $a_1=a$, $a_q=b$, set $x_1=x$, $y_q=y$, and for every $i \in \{1,\dots,q-1\}$ we choose edges $e_i:=y_1x_{i+1} \in E[V'_{a_i},V'_{a_{i+1}}]$ of alternating color starting with blue.
		Note that $y_1 \neq x$ as $y_1x_2$ is a blue edge, and therefore $y_1$ is blue.
		Similarly, $y\neq x_q$ as $y_{q-1}x_q$ is a red edge, and therefore $x_q$ is red.

		By Hamilton laceability, we obtain that for every $i \in \{0, \dots, q-1\}$ there is a Hamilton $x_iy_i$-path $P^i$.
		We conclude by noting that the concatenation $P^1P^2\cdots P^{q-1}P^{q}$ is a Hamilton $xy$-path.
		
		\noindent\textbf{Case 2:} $q \equiv 1 \pmod{4}$.
		By Lemma~\ref{lem:alternating_coloring}, there is an alternating Hamilton $ab$-path $a=a_1,\dots,a_q=b$ in~$\bar{K}_q$ that starts with a blue edge.
		As $q'$ is even in this case, it is easy to see that the edges between $V_c$ and $V_d$ can only be of one color.
		Moreover, by the coloring scheme (*), this color is the same as the color of the edge $cd$ in $\bar{K}_q$. 
		All of the above implies that we can choose edges $e_i:=x_iy_{i+1} \in E[V'_{a_i},V'_{a_{i+1}}]$ of alternating colors starting with blue, for $i \in \{1,\dots,q-1\}$.
		We additionally set $x_1 = x$ and $y_q = y$.
		As before, since $y_1x_2$ is a blue edge, $y_1$ is blue, and hence $y_1 \neq x$.
		Similarly, since $y_{q-1}x_q$ is a red edge, $x_q$ is red, and hence $y\neq x_q$.

		By Hamilton laceability, we obtain that for every $i \in \{0, \dots, q-1\}$ there is a Hamilton $x_iy_i$-path $P^i$.
		We again conclude by noting that the concatenation $P^1P^2\cdots P^{q-1}P^{q}$ is a Hamilton $xy$-path.
		\begin{figure}[h]
			\centering
			\includegraphics[width=\textwidth,page=1]{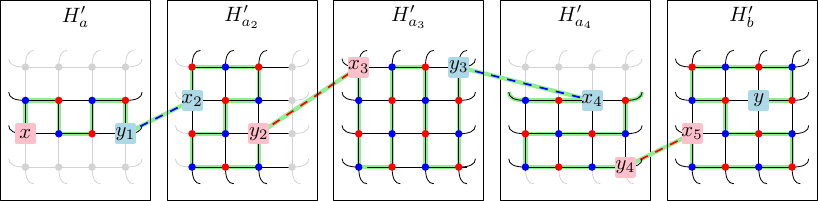}
			\caption{Proof illustration of Lemma~\ref{lem:order_components}.}
			\label{fig:subcase1}
		\end{figure}
		\end{proof}

We are now ready to prove that $H$ is Hamilton connected.

\begin{lem}
	\label{lem:H_connected}
	For every $q\geq 3$ the graph $H(q)$ is Hamilton connected.
\end{lem}
\begin{proof}
	We consider three cases.
	
	\begin{figure}[ht]
	\centering
	\includegraphics[width=.5\textwidth]{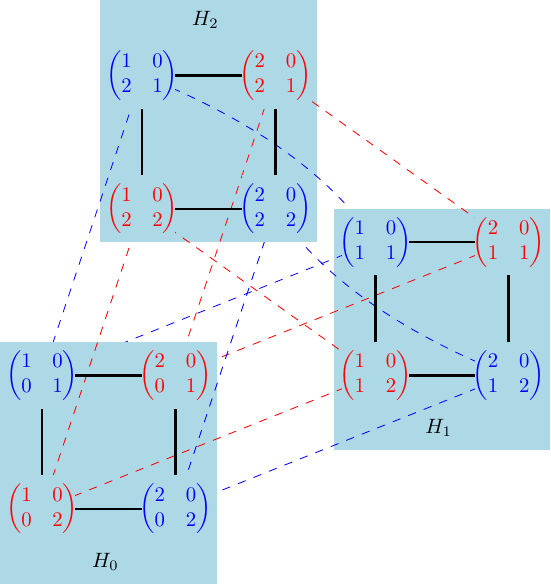}
	\caption{The graph~$H$ for $q=3$ and $n=2$.
	For these particular values of~$q$ and $n$ the graph~$H$ is isomorphic to $C_4 \boxprod C_3$.}
	\label{fig:Hgraph}
	\end{figure}
	
	\noindent\textbf{Case 1:} \textit{$q=3$.}
		Here, $H \cong C_4 \boxprod C_3$ (see Fig.~\ref{fig:Hgraph}) which is Hamilton connected by Lemma~\ref{lem:cycleproduct}.

	\noindent\textbf{Case 2:} \textit{$q=2^k$ for $k\geq 2$.}
		We apply Lemma~\ref{lem:join} on $H$ with the partition $\{V_a \, : \,  a \in \FF_q \}$.

		To see that (1) holds, note that for every $a \in \FF_q$ we have that $H_a \cong C_{2^k-1} \boxprod C_{2^k-1}$, and by Lemma~\ref{lem:cycleproduct} we conclude that $H_a$ is Hamilton connected.

		For every $a\in \FF_q$ and $x \in V_a$ we have that $x=\begin{pmatrix}\alpha^i & 0\\ \alpha^j a & \alpha^j \end{pmatrix}$ for some $i,j \in \{0,\dots, q-2\}$.
		Thus, when adding the first row to the second we obtain a matrix in $V_{a+\alpha^{j-i}}$.
		Since $\alpha^{j-i} \neq 0$, we conclude that (2) holds.

		Finally, for every $a,b\in \FF_q$ we have $q-1\geq 3$ disjoint edges between $V_a$ and $V_b$, thus (3) also holds.
		Therefore, applying Lemma~\ref{lem:join} we conclude that $H$ is Hamilton connected.

	\noindent\textbf{Case 3:} \textit{$q\geq 5$ and odd.}
	In this case for every $a \in \FF_q$ the graph $H_a \cong C_{q-1} \boxprod C_{q-1}$ is bipartite with partitions given by the colors red and blue.
	Thus, $H_a$ is Hamilton laceable by Lemma~\ref{lem:cartesian}.

	Let $x \in V_a$ and $y \in V_b$. 
	We aim to show that there is a Hamilton $xy$-path in $H$.
	Let $i_x,i_y,j_x,j_y \in \{0,\dots, q-2\}$ and $a,b \in \FF_q$ such that
	$x=\begin{pmatrix} \alpha^{i_x} & 0 \\ \alpha^{j_x} a & \alpha^{j_x}\end{pmatrix}$ and $y=\begin{pmatrix}
			\alpha^{i_y} & 0 \\ \alpha^{j_y} b & \alpha^{j_y}\end{pmatrix}$.
	Note that either $i_x \neq i_y$ or $j_x \neq j_y$.
	We assume that $i_x \neq i_y$ as the case $j_x\neq j_y$ is analogous.
	Furthermore, we assume without loss of generality that $\col(x)$ is red.

	\begin{itemize}
		\item \textbf{Subcase 3.1:} $a\neq b$ and $\col(x)\neq \col(y)$.
		This case is a direct consequence of Lemma~\ref{lem:order_components} by using that $H$ is a structured induced subgraph of $H$.
	
		\item \textbf{Subcase 3.2:} $a\neq b$ and $\col(x)=\col(y)$.
		We consider the following set
		\[V^* := \left\{ \begin{pmatrix} \alpha^{i} & 0 \\ \alpha^{j}a & \alpha^{j}  \end{pmatrix}  : i \in \{0,\dots, q-2\}, j \in \{j_x,\dots,j_x+ (q-3)/2\} \right\} \subseteq V_a \]
		and note that $x \in V^*$.
		Choose a blue vertex $y'\in V^*$ and a vertex $x' \in V(H) \setminus (V_a\cup V_b)$ such that $y'x' \in E[V_a,V(H)\setminus (V_a\cup V_b)]$.
		Note that $H^*:=H[V^*] \cong C_{q-1} \boxprod P_{\frac{q-1}{2}}$ and thus there exists $P^1$ a Hamilton $xy'$-path in $H^*$.
		Additionally, since $y'$ is blue, then $x'$ is also blue.
		It is easy to check that $H\setminus V^*$ is structured. As a consequence, Lemma~\ref{lem:order_components} implies that there is a Hamilton $x'y$-path $P^2$ in $H\setminus V^*$.
		We conclude by noting that the concatenation $P^1P^2$ is a Hamilton $xy$-path in $H$.

		\begin{figure}
			\centering
			\begin{subfigure}{.45\textwidth}
			\centering
			\includegraphics[page=1]{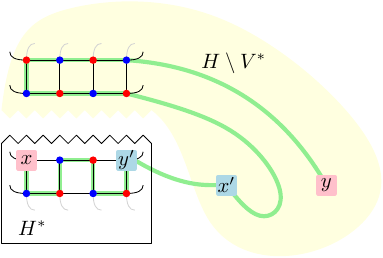}
			\caption{}
			\end{subfigure}\hfill
			\begin{subfigure}{.45\textwidth}
				\centering
				\includegraphics[page=2]{basecase.pdf}
				\caption{}
				\end{subfigure}
			
			\begin{subfigure}{.9\textwidth}
				\centering
				\includegraphics[page=3]{basecase.pdf}
				\caption{}
			\end{subfigure}
			\caption{Illustrations for Lemma~\ref{lem:H_connected}'s proof in subcase (a) 3.2, (b) 3.3, and (c) 3.4.}
		\end{figure}

		\item \textbf{Subcase 3.3:} $a=b$ and $\col(x)=\col(y)$.
		Choose a partition of $\{0,\dots, q-2\}$ into sets $I_x,I_y$ such that each set is a (cyclic) interval, $i_x \in I_x$, $i_y \in I_y$, and $|I_x|=|I_y|=\frac{q-1}{2}$.
		We now define
		\[V^* := \left\{ \begin{pmatrix} \alpha^{i} & 0 \\ \alpha^{j}a & \alpha^{j}  \end{pmatrix}  : i \in I_x, j \in \{0,\dots,q-2\} \right\} \subseteq V_a,\]
		define $H^* := H[V^*]$.
		Note that $y \notin V^*$, thus we can continue the proof exactly as in subcase 3.2.
		\item \textbf{Subcase 3.4:} $a=b$ and $\col(x) \neq \col(y)$.
		As before, we choose a partition of $\{0,\dots, q-2\}$ into sets $I_x,I_y$ such that each set is a (cyclic) interval, $i_x \in I_x$, $i_y \in I_y$, and $|I_x|=|I_y|=\frac{q-1}{2}$, let
		\[V^1 := \left\{ \begin{pmatrix} \alpha^{i} & 0 \\ \alpha^{j}a & \alpha^{j}  \end{pmatrix}  : i \in I_x, j \in \{0,\dots,q-2\} \right\} \subseteq V_a,\]
		and let $H^1 := H[V^1]$.
		Choose a blue vertex $y_1\in V^1$ and a vertex $x_2 \in V(H) \setminus V_a$ such that $y_1x_2 \in E[V_a,V(H)\setminus V_a]$.
		Let $i^*,j^* \in \{0,\dots, q-2\}$, $c \in \FF_q$ such that $x_2 = \begin{pmatrix}
			\alpha^{i^*} & 0 \\
			\alpha^{j^*}c & \alpha^{j^*}.
		\end{pmatrix}$.
		We define
		\[V^2 := \left\{ \begin{pmatrix} \alpha^{i} & 0 \\ \alpha^{j}c & \alpha^{j}  \end{pmatrix}  : i \in \{0,\dots,q-2\}, j \in \{j^*,\dots,j^*+(q-3)/2\} \right\} \subseteq V_c,\]
		and let $H^2 := H[V^2]$.
		Note that $H^1 \cong H^2 \cong C_{q-1} \boxprod P_{\frac{q-1}{2}}$.
		Thus, by Lemma~\ref{lem:pathxcycle}, there exists $P^1$ a Hamilton $xy_1$-path in $H^1$ and $P^2$ a Hamilton $x_2y_2$-path in $H^2$.
		
		We now show that $H \setminus (V^1 \cup V^2)$ is structured. 
		Let $k \in \{0, \dots, q-2\}$ be one of the two common points of $I_x$ and $I_y$.
		Further suppose that $a - c = \alpha^z$ for some $z \in \{0, \dots, q-2\}$.
		Since $(I_x, I_y)$ is a partition of $\{0, \dots, q-2\}$, and since the two sets have equal size, either $k + z$ or $k - z$ is in $I_x$.
		Without loss of generality, suppose this is $k+z$.
		Then $\begin{pmatrix} \alpha^{k+z} & 0 \\ \alpha^{k}a & \alpha^{k}  \end{pmatrix} \begin{pmatrix} \alpha^{k+z} & 0 \\ \alpha^{k}c & \alpha^{k}  \end{pmatrix}$ is an edge between $V_a \setminus V^1$ and $V_c \setminus V^2$, by (p3).
		The other conditions of Definition~\ref{def:structured} are easy to check.
		Hence, $H \setminus (V^1 \cup V^2)$ is structured.
		
		As a consequence Lemma~\ref{lem:order_components}, the above implies that there is a Hamilton $x_2y$-path $P^3$ in $H\setminus (V^1\cup V^2)$.
		We conclude by noting that the concatenation $P^1P^2P^3$ is a Hamilton $xy$-path in $H$.
		\qedhere
	\end{itemize}
\end{proof}

We are now ready to prove Lemma~\ref{lem:base_case}.

\begin{proof}[Proof of Lemma~\ref{lem:base_case}]	
	Let us denote~$S = \left\{ (a, 1) \,:\, a\in\FF_q \right\} \cup \left\{(1, 0)\right\} \subseteq \FF_q^2$.
	For any nonzero vector~$u\in S$, let us define $V^u \subseteq V(G)$ to be the set of all matrices in~$G$ with the first row in $\vspan{u}$.
	Note that each $V^u$ corresponds to a unique component of~$G'$, with $V^{(1, 0)} = V(H)$.
	
	We apply Lemma~\ref{lem:join} with partitioning $\left\{V_u \,:\, u\in S\right\}$.
	
	We begin by simplifying our arguments using symmetry.
	Note that any two components of~$G'$ are isomorphic via $f_A: x \mapsto xA$ for some~$A\in GL(2, q)$.
	Moreover, for any~$A$ the mapping~$f_A$ is an automorphism of~$G$ that preserves operations on the edges; i.e., if an edge corresponds to the operation~$M_1$, then it will be mapped to some edge that also corresponds to~$M_1$.
	In particular, $H$ is isomorphic to every other component.
	
	For the condition~(1) of Lemma~\ref{lem:join}, we know that~$H$ is Hamilton connected by Lemma~\ref{lem:H_connected}, so by isomorphism the same holds for every~$G[V^u]$.	
	The condition~(2) is satisfied simply by adding the second row to the first row.
	For the condition~(3), we first observe that if there is an edge between~$G[V^u]$ and~$G[V^{u'}]$,	there are actually at least~$q-1$ disjoint edges as we can multiply both matrices by~$\alpha^i$.	
	By isomorphism, it is enough to show that there is an edge between~$H$ and any~$V^u$ distinct from~$V^{(1, 0)}$.
	Since~$u = (a, 1)$ for some~$a \in \FF_q$, we can use the edges between $\begin{pmatrix} 1 & 0 \\ a-1  & 1 \end{pmatrix}\in V(H)$ and
	$\begin{pmatrix} a & 1 \\ a-1  & 1 \end{pmatrix}\in V^{(a, 1)}$.
\end{proof}

\section{Alternating path of two-edge-colored complete graph}
\label{sec:alternating_coloring}

In this section, we prove Lemma~\ref{lem:alternating_coloring}, which is needed in the proof of Lemma~\ref{lem:order_components}.

We start with a brief recap of the context needed for this lemma.
Suppose $q$ is a prime power and $q \equiv 1$ (mod 4).
We remind the reader that $\alpha$ is a generator of the multiplicative group of $\mathbb{F}_q$.
Recall that $\bar{K}_{q}$ is the complete graph on the vertex set $\mathbb{F}_q$ with edges colored by the following scheme:
\begin{itemize}
	\item[(*)] For $x$ and $y$ in $\mathbb{F}_q$, the edge $xy$ has color red (blue), if there exists an odd (even)~$z \in \mathbb{Z}$, such that $x - y = \alpha^z$.
\end{itemize}
Our goal is to find an alternating Hamilton path between two prescribed vertices~$a$ and~$b$ with a prescribed color of the edge incident to $a$.

We begin by arguing that $\bar{K}_{q}$ is well-defined.
Let 0 be the additive identity and 1 be the multiplicative identity of $\mathbb{F}_q$.
By the definition of $\alpha$, the nonzero elements of $\mathbb{F}_q$ are exactly $\alpha^0, \dots, \alpha^{q-2}$.
Furthermore, $\alpha^{i} = \alpha^{q-1+i}$ for all integers $i\in\mathbb{Z}$.
Since $q$ is odd, we conclude that if $\alpha^z = \alpha^{z'}$ for some $z, z'$, then $z$ and $z'$ have the same parity.
Thus, for $x$ and $y$ in $\mathbb{F}_q$, there exists a unique $p \in \{0,1\}$, such that if $x - y = \alpha^z$ then $z \equiv p$ (mod 2).
Further, since $\alpha^{(q-1)/2} = -1$, if $x-y = \alpha^z$, then $y-x = \alpha^{z'}$ for $z' = z + (q-1)/2$.
As $q \equiv 1$ (mod 4), $z$ and $z'$ have the same parity.
Therefore, the color of each edge of $\bar{K}_{q}$ is well-defined.

The problem of finding an alternating cycles/paths in a graph has a long history and a wide range of applications; see a survey by Bang-Jensen and Gutin~\cite{MR1439259}.
The earliest result on alternating trail in 1968 by B\'{a}nkfalvi and B\'{a}nkfalvi~\cite{MR0233731} gave a characterization of the two-edge-colored complete graphs that have an alternating Hamilton cycle.
However, these graphs must have an even number of vertices, and as such, we cannot readily apply the result.
Specific to alternating Hamilton paths, Bang-Jenson, Gutin, and Yeo~\cite{MR1609957} showed the following necessary and sufficient condition.
A two-edge-colored complete graph $G$ has an alternating Hamilton $ab$-path, if and only if $V(G)$ can be partitioned into disjoint subsets $(V_0, V_1, \dots, V_t)$, for some $t \geq 0$, such that there are an alternating $ab$-path that spans $V_0$ and an alternating cycle that spans each of $V_i$ for $i \in [t]$.
However, it is not evident how we can apply the result above in our setting.
Furthermore, we also specify the color of the first edge of the path, which is not guaranteed by the statement above.
Therefore, we provide a direct and constructive proof of an alternating Hamilton path in our special complete graph.

In the following proof, we use the observation that the operation of adding a constant to all vertex labels preserves edge colors, since the difference between any two vertices remains the same.

\alternatingcoloring*

\begin{proof}
    For $i \in \{0, \dots, q-1\}$, define $v_i := \sum^{i}_{j=0} \alpha^{j}$.
	Note that $v_{q-2} = 0$, and $v_0 = v_{q-1} = 1$.
    It is easy to see that $(v_0, \dots, v_{q-2})$ forms an alternating cycle~$C$ in~$\bar{K}_{q}$.
    The only missing vertex in~$C$ is $u:=-(\alpha-1)^{-1}$.
    Indeed, if this vertex is on the cycle, then for some $t$, we must have $(\alpha^{t+1}-1)(\alpha-1)^{-1} = -(\alpha-1)^{-1}$, which implies $\alpha^{t+1} = 0$, a contradiction with the fact that $\alpha$ generates nonzero elements of $\mathbb{F}_q$.

    For any $i \in \{0, \dots, q-2\}$ we have
    \[
        v_i - u = (\alpha^{i+1}-1)(\alpha-1)^{-1} + (\alpha-1)^{-1} = \alpha^{i+1} (\alpha-1)^{-1}.
    \]
    By a similar argument, we have that
    \[
        v_{i+1} - u = \alpha^{i+2} (\alpha-1)^{-1} = \alpha(v_i - u).
    \]
    Thus, by the coloring scheme (*), $uv_i$ and $uv_{i+1}$ have different colors.

\begin{claim*}
    For any vertex~$v$ in $C$, there exists an alternating Hamilton $uv$-path such that on the path, $u$ is incident to an edge whose color is different from that of $uv$.
\end{claim*}
{\makeatletter
\let\par\@@par
\par\parshape0
\everypar{}
\begin{wrapfigure}{r}{0.5\textwidth}
   \centering
    \includegraphics[scale=0.9]{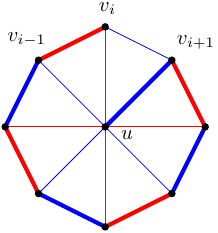}
    \captionof{figure}{Illustration of the claim's proof. The outer cycle is the cycle~$C$, and the bold edges indicate an alternating Hamilton path.}
    \label{fig:wheel}
  \end{wrapfigure}
\vspace{-1.4em}\noindent
\begin{claimproof}
    Suppose $v = v_i$ for some $i \in \{0, \dots, q-2\}$.
    By the argument above, $uv_{i-1}$ and $uv_{i+1}$ have the same color.
    Further, since $C$ is an alternating cycle, $v_{i-1}v_i$ and $v_iv_{i+1}$ have different colors.
    Hence, one of these two edges have the same color as $uv_{i-1}$ and $uv_{i+1}$.
    Without loss of generality, suppose this edge is $v_iv_{i+1}$.
    Then we have $(u, v_{i+1}, v_{i+2} \dots, v_{q-2}, v_0, \dots, v_{i-1}, v_i)$ is the desired alternating Hamilton path.
    See Figure~\ref{fig:wheel} for an illustration. 
\end{claimproof}

}%


    Consider adding $b-u$ to all vertex labels.
	The missing vertex from the cycle~$C$ above is now $b$.
    By the claim above, we obtain an alternating Hamilton $ba$-path, such that the incident edge to $b$ has different color than that of $ba$.
    Since $q$ is odd, this implies that the edge incident to $a$ on this path has the same color as $ba$.
    Next, we add $a-u$ to all original vertex labels.
	The missing vertex from~$C$ is now $a$.
    Again by the claim above, we obtain another alternating Hamilton $ab$-path, such that the incident edge to $a$ has different color than that of $ab$.

    Since the two alternating Hamilton paths above have different colors for the edge incident to $a$, the lemma follows.
\end{proof}

\section{Algorithmization}

\label{sec:algorithmic}

The proof of Theorem~\ref{thm:ham_connected} can be easily turned into an algorithm that computes a Hamilton path in $G(n,q,T)$ running in time polynomial in $|GL(n,q)|$.
This can be obtained by a straightforward recursion based on the joining lemma.
More specifcally, the main idea of the proof is to split the graph and proceed recursively.
Close examination of the proof of the joining lemma and its applications along our proof shows that such a recursion can be computed in time polynomial in $|GL(n,q)|$; we omit the details.

However, the typical goal from a generation perspective is to have an algorithm that outputs objects one by one with a small delay and preprocessing time.
Here, the \defi{delay} is the worst-case time the algorithm takes between consecutively generated objects and the \defi{preprocessing time} is the time before generating any objects.
Thus, the natural objective from generation perspective for invertible matrices is an enumeration algorithm running in delay $\poly(n,q) := n^{\cO(1)}q^{\cO(1)}$ with $\poly(n,q)$ preprocessing.
Note that such an algorithm immediately gives a solution to computing Hamilton paths in $G(n,q,T)$ in time $\poly(n,q)|GL(n,q)|$.

The naive implementation of our inductive proof uses a call stack that needs space exponential in $n$ and takes exponential time in $n$ to put the recursive calls in the stack.
Despite that, it is still possible to obtain a polynomial delay algorithm by following the recursive structure of the main proof.
As highlighted in the stack approach, we cannot store \emph{all} the information given by the recursion.
Instead, we only store information related to the \emph{current path} in the recursion tree.
More specifically, if we are at a vertex $z$, we trace back the $\ell$ recursive calls, each utilizing the joining lemma.
The $i$-th call indicates a pair of a source $x_i$ and a target $y_i$ for which we traverse a Hamilton path.
For every $i \in [\ell]$ we store $x_i$, $y_i$ and a small amount of extra bits serving as a compressed history.
It turns out that this is enough information to reconstruct the path of $z$ in the recursion tree and decide how to proceed; more details are given in Appendix~\ref{app:algo}.

\section{Open questions}
We conclude with several remarks and open questions.
\begin{enumerate}
\item {\bf Algorithmization.}
The proof of Theorem~\ref{thm:ham_connected} heavily relies on an induction that can be called multiple times in a single step and requires keeping history in the memory.
Thus, it is not suitable for efficient algorithm.
Is there a generating algorithm that achieves $\cO(n)$ delay? Is there a simple greedy algorithm?
\item {\bf Non-bypass transition graphs.}
The transition graph specifies which row additions and subtractions we allow.
For our induction step we need that it is a bypass transition graph.
Does Theorem~\ref{thm:ham_connected} hold for any strongly connected transition graph, in particular for the directed $n$-cycle?
We verified by computer that the result holds for the directed cycle if $q=2$ and $n=3$.
\item {\bf Other generators.}
We used multiplication/division by $\alpha$ as an elementary operation for any row.
Can we adapt our methods for multiplication/division of only a fixed row?
On the other hand, we may allow general operations and consider any generator of the group $GL(n,q)$. In particular the generator $\{M_2 A_{1,n},P_{(2,\dots,n,1)}\}$ of size $2$, where $P_{(2,\dots,n,1)}$ refers to the permutation matrix corresponding to the permutation $2\dots n1$~\cite{MR1007258}.
\item {\bf Subgroups of $GL(n,q)$.} As an intermediate step, we show that Cayley graphs of certain subgroups of $GL(n,q)$ are Hamilton connected.
Can we prove it for other subgroups that correspond to given restrictions of matrices?
\item {\bf Symmetric Hamilton cycles.}
Instead of Hamilton connectivity we may ask for Hamilton cycles that are preserved under a large cyclic subgroup of automorphisms.
This problem was recently studied by Gregor, Merino, and Mütze~\cite{gregor_et_al:LIPIcs.MFCS.2022.54} for several highly symmetric graphs.
The graphs considered here are also highly symmetric.
For example, we know that for $q=2$, $n=3$, and the complete transition tree there exists a 24-symmetric Hamilton cycle. It can be shown that the automorphism group $\Aut(\Cay(GL(n,2),\{A_{ij} : 1\le i<j \le n\}))$ is isomorphic to $S_n \times (\mathbb{Z}_2 \ltimes GL(n,2))$ for every $n\ge 3$ and the dihedral group $D_{12}$ for $n=2$, so there can be highly symmetric Hamilton cycles.
\item {\bf Matrices over rings.}
Another natural extension is to explore if our results extend to invertible matrices in the ring setting.
This is particularly interesting for cyclic rings; i.e., checking Hamiltonicity of Cayley graphs of invertible matrices in $\ZZ_k$ for $k \in \NN$.
Naturally, the methods will highly depend on the chosen generators for which there does not seem to be an obvious choice.
\item {\bf Alternating Hamilton paths in 2-colored $K_{2n+1}$.}
Despite our efforts and many existing results on properly colored Hamilton cycles in complete graphs (see a survey \cite{MR1439259}), we did not find an answer to the following question.
Is it true that the complete graph $K_{2n+1}$ with $2$-colored edges so that every vertex is incident with exactly $n$ edges of each color contains an alternating Hamilton path between any two vertices?
\end{enumerate}

\bibliographystyle{plain}
\bibliography{refs}

\section{Base case checking}

\label{app:sage}

We provide the following SageMath code to check the base cases in Lemma~\ref{lem:base_case32}.

\definecolor{codegreen}{rgb}{0,0.6,0}
\definecolor{codegray}{rgb}{0.5,0.5,0.5}
\definecolor{codepurple}{rgb}{0.58,0,0.82}
\definecolor{backcolour}{rgb}{0.95,0.95,0.92}

\lstdefinestyle{mystyle}{
    backgroundcolor=\color{backcolour},
    commentstyle=\color{codegreen},
    keywordstyle=\color{magenta},
    numberstyle=\tiny\color{codegray},
    stringstyle=\color{codepurple},
    basicstyle=\ttfamily\footnotesize,
    breakatwhitespace=false,
    breaklines=true,
    captionpos=b,
    keepspaces=true,
    numbers=left,
    numbersep=5pt,
    showspaces=false,
    showstringspaces=false,
    showtabs=false,
    tabsize=2,
	indent=2ex
}

\begin{lstlisting}[language=Python, tabsize=2,backgroundcolor=\color{backcolour},
    commentstyle=\color{codegreen},
    keywordstyle=\color{magenta},
    numberstyle=\tiny\color{codegray},
    stringstyle=\color{codepurple},basicstyle=\ttfamily\footnotesize]
def base_case_graph(T):
	# Builds the base case graph G(3,2,T)
	GL_group = GL(3,GF(2))
	V = GL_group.list()
	elem= []
	for i in range(3):
		for j in range(3):
			if (i,j) in T:
				elem.append(elementary_matrix(GF(2),3,row1=i, row2=j, scale=1))
	G = Graph()
	for A in V:
		for E in elem:
			if str(A)!=str(E*A):
				G.add_edge(str(A),str(E*A),label=str(E))
	return G

def check_base_cases():
	# Checks all base cases
	bypass_graphs = [
		[(0,1),(1,0),(1,2),(2,1)], [(0,1),(1,0),(1,2),(2,0)],
		[(0,1),(1,0),(2,1),(0,2)], [(0,1),(1,0),(2,0),(0,2)]
	] # All bypass graphs on 3 nodes
	I = str(matrix.identity(3))
	# Only check ham paths from identity due to vertex-transitivity
	for T in bypass_graphs:
		G = base_case_graph(T)
		for v in G.vertices():
			if v!=I:
				if not G.hamiltonian_path(I,v):
					return False
	return True

print(check_base_cases())
\end{lstlisting}

\section{Further details on algorithmization}
\label{app:algo}
In this appendix, we elaborate further our idea on how we can turn the proof of Theorem~\ref{thm:ham_connected} to a poly($n,q$)-delay generation algorithm.
In particular, we give an example of how we can implement an application of the joining lemma (Lemma~\ref{lem:join}) in the proof and discuss the preprocessing time.

We start with the preprocessing algorithm, which essentially handles the base cases.
Specifically, given $q$, we turn the proof of either Lemma~\ref{lem:base_case32} (if $q = 2$) or Lemma~\ref{lem:base_case} (if $q > 2$) into a preprocessing algorithm to produce Hamilton paths for all pairs of starting and ending invertible matrices in the relevant base case.
The number of invertible matrices in this base case is $\cO(q^9)$.
As discussed in Section~\ref{sec:algorithmic}, our algorithm takes time polynomial to this number, and hence, the runtime is polynomial in $q$.
Moreover, it takes $\cO(q^{18})$ space to store all the Hamilton paths for all pairs of starting and ending matrices.
With this preprocessing, for every matrix in the base case, we can query the next matrix in the relevant Hamilton path in constant time.

\begin{figure}
	\centering
	\includegraphics[width=\textwidth]{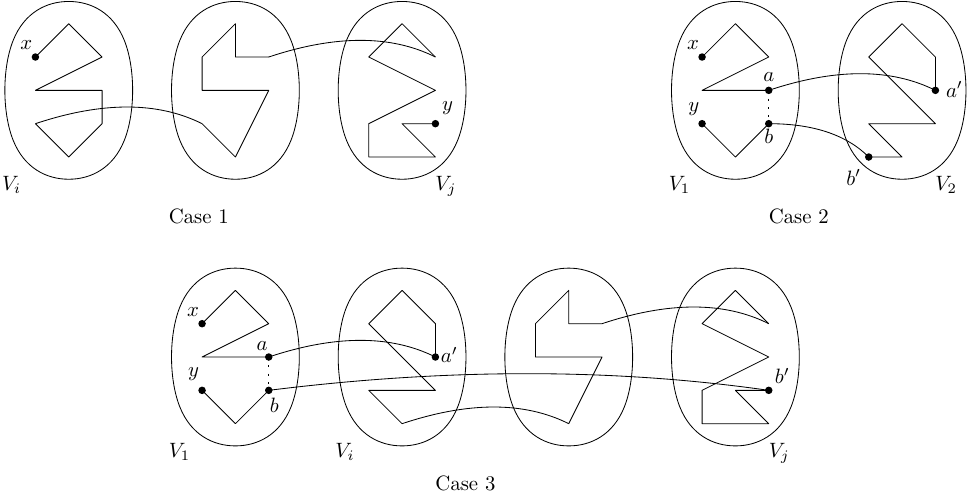}
	\caption{
		Illustration of the proof of Lemma~\ref{lem:join}.
	}
	\label{fig:join}
\end{figure}

Next, we recap briefly the proof idea of Lemma~\ref{lem:join}.
Suppose we want to find a Hamilton $xy$-path in $G = V_1 \cup \dots \cup V_k$ such that the three conditions stated in the lemma hold.
In the proof, we essentially consider three cases. (See Figure~\ref{fig:join} for an illustration.)

\noindent
\textbf{Case 1}: \textit{$x$ and $y$ are in the different sets $V_i$ and $V_j$, respectively}.
We consider an ordering of the subsets starting with $V_i$ and ending with $V_j$.
Then the desired Hamilton path visits all vertices in each subset before moving to the next subset in this ordering.

\noindent
\textbf{Case 2}: \textit{$x$ and $y$ are in the same set $V_1$, and $k=2$}.
In this case, we traverse along the Hamilton $xy$-path in $V_1$ until we encounter an edge $ab$ such that the neighbors $a'$ and $b'$ of $a$ and $b$ in $V_2$ , respectively, are distinct.
We then traverse the Hamilton $a'b'$-path in $V_2$ before coming back to the Hamilton $xy$-path in $V_1$.

\noindent
\textbf{Case 3}: \textit{$x$ and $y$ are in the same set $V_1$, and $k>2$}.
Here, we traverse along the Hamilton $xy$-path in $V_1$ until we encounter an edge $ab$ such that the neighbors of $a'$ and $b'$ are in some subsets $V_i$ and $V_j$.
We use the same technique in Case 1, where we consider an ordering of the subsets other than $V_1$ such that this ordering starts with $V_i$ and ends with $V_j$.
Then we traverse all vertices in these subsets in this ordering, before coming back to the Hamilton $xy$-path in $V_1$.

With this recap in mind, we consider the inductive step presented in Section~\ref{sec:main_proof}.
The high level idea of the proof for this step is to apply Lemma~\ref{lem:join} to join the different subgraphs $G_n[(S_u, *)]$ for all one-dimensional subspaces $\langle u \rangle$.
We can package this as a function $\textsc{Join}_{\text{F}}(n,q,T,x,y,z,\alpha)$ that takes as input a natural number $n$, a prime power $q$, a bypass transition graph $T$ on $n$ vertices, three invertible $n \times n$ matrices $x,y,z$, and a binary number $\alpha$.
(The subscript F refers to the application of the joining lemma to the \emph{full} graph $G(n,q,T)$.)
Then it returns the next matrix after the current matrix $z$ in a Hamilton $xy$-path of $G_n:=G(n,q,T)$.
The binary number $\alpha$ is an auxiliary variable that takes the value 1 if and only if $x$ and $y$ are in the same subgraph $G_n[(S_{\bar{u}}, *)]$ for some $\bar{u}$, and before visiting $z$, we have visited all matrices in $G_n \setminus G_n[(S_{\bar{u}}, *)]$.
This variable is a technical detail to accommodate the fact that we come back to the subset containing $x$ and $y$ after traversing all other subsets in Cases 2 and 3 in the proof of the joining lemma.

Based on the conditions in the statement and the proof of the joining lemma, we have the following subroutines.
\begin{itemize}
	\item For a matrix $x$ in $G_n$, the subroutine $\textsc{idx}(n, q, x)$ returns a one-dimensional vector $u$ such that $x \in (S_u,*)$.
	\item (First condition of Lemma~\ref{lem:join})
		For a vector $u \in F^n_q$, three $n \times n$ matrices $x, y, z \in (S_u,*)$, the subroutine $\textsc{Join}_{\text{EL}}(n,q,T,u,x,y,z)$ returns the next matrix after $z$ in a Hamilton $xy$-path of $G_n[(S_u,*)]$.
		(Here, the subscript EL refers to the application of the joining lemma to $G_n[(S_u,*)]$, that is, the induced subgraph on the set of matrices with certain constraints on all the rows \emph{except the last one}.)
	\item (Second condition of Lemma~\ref{lem:join})
		For a matrix~$x$ in $G_n$, the subroutine $\textsc{OtherSet}(n,q,T,x)$ returns the matrix $y$ such that $xy$ is an edge in $G_n$ and $\textsc{idx}(n, q, x) \neq \textsc{idx}(n, q, y)$.
	\item (Third condition of Lemma~\ref{lem:join})
		For two vectors $u, u' \in F^n_q$, the subroutine $\textsc{Conn}(n,q,T,u,u')$ returns three disjoint edges $x_1y_1$, $x_2y_2$, and $x_3y_3$ such that $\textsc{idx}(n,q,x_i) = u$ and $\textsc{idx}(n,q,x_j) = u'$ for $i \in [3]$.
	\item For two matrices $x, y$ such that $\textsc{idx}(n, q, x) \neq \textsc{idx}(n, q, y)$ and a vector $u \in F^n_q$, the subroutine $\textsc{SubsetOrderDiff}(n,q,x,y,u)$ returns the next subset after $(S_u,*)$ in the ordering of the subsets $(S_w,*)$ for all one-dimensional subspaces $\langle w \rangle$ such that the ordering starts with $(S_{\textsc{idx}(x)},*)$ and ends with $(S_{\textsc{idx}(y)},*)$.
	\item For a vector $u \in F^n_q$ and four matrices $x, y, a', b'$ such that $\textsc{idx}(n, q, x) = \textsc{idx}(n, q, y) = u'$ for some $u$ and $\textsc{idx}(a'), \textsc{idx}(b') \neq u'$, the subroutine $\textsc{SubsetOrderSame}(n,q,x,y,a',b',u)$ returns the next subset after $(S_u,*)$ in the ordering of the subsets other than $(S_{u'},*)$ such that the ordering starts with $(S_{\textsc{idx}(a')},*)$ and ends with $(S_{\textsc{idx}(b')},*)$.
\end{itemize}

With all these subroutines, it is easy to complete $\textsc{Join}_{\text{F}}(n,q,T,x,y,z,\alpha)$ by following the proof of Lemma~\ref{lem:join}.
It remains to argue about the running time of these subroutines.

The subroutine $\textsc{idx}$ can be done by solving a system of linear equations and hence in polynomial time.
The subroutines $\textsc{OtherSet}$ and $\textsc{Conn}$ can be implemented straightforwardly from our proof in Section~\ref{sec:main_proof}, where we give explicit constructions of the outputs.

For $\textsc{SubsetOrderDiff}$ and $\textsc{SubsetOrderSame}$, we first define a canonical ordering of all subsets $(S_w,*)$ for all one-dimensional subspaces $\langle w \rangle$.
For each such subspace, we choose its representative to be the unique vector whose first nonzero element is 1.
Then an ordering of the subsets can be described by an ordering of all vectors whose first nonzero element is 1.
We choose the canonical ordering to be the lexicographical ordering of these vectors.
It is then easy to see that in polynomial time, we can obtain the next vector in the ordering.
Now, the ordering for $\textsc{SubsetOrderDiff}(n,q,x,y,u)$ is obtained from the canonical ordering above by moving $(S_{\textsc{idx}(x)},*)$ to the beginning and $(S_{\textsc{idx}(y)},*)$ to the end.
In other words, we iteratively check the subset after $(S_u,*)$ in the canonical ordering until we obtain a subset that is neither $(S_{\textsc{idx}(x)},*)$ nor $(S_{\textsc{idx}(y)},*)$.
The ordering for $\textsc{SubsetOrderSame}(n,q,x,y,a',b',u)$ can be obtained in a similar fashion.

The remaining subroutine is $\textsc{Join}_{\text{EL}}(n,q,T,u,x,y,z)$.
Here, we can implement this subroutine by using a similar approach as how we implement $\textsc{Join}_{\text{F}}(n,q,T,x,y,z,\alpha)$.
Note that similarly to $\textsc{Join}_{\text{F}}$, we may need an additional binary variable in the input of $\textsc{Join}_{\text{EL}}$ to accommodate a technical detail in the proof of Lemma~\ref{lem:join}.

\end{document}